\documentclass[10pt,amssymb]{article}

\usepackage{eucal}
\usepackage{amsfonts,amsmath,amssymb,latexsym}

\begin{document}

\newtheorem{theorem}{Theorem}
\newtheorem{prop}[theorem]{Proposition}
\newtheorem{lemma}[theorem]{Lemma}
\newtheorem{corollary}[theorem]{Corollary}
\newtheorem{definition}{Definition}[section]
\newtheorem{proof}{proof}

\def\Ss{\mathbb{S}}
\def\Mm{\mathbb{M}^2({\kappa})}
\def\Cc{\mathcal{C}}
\def\Rr{\mathbb{R}}
\def\Pp{\mathbb{P}_t}
\def\dd{\mathrm{d}}
\def\mt{\mathcal}
\def\ae{\langle}
\def\ad{\rangle}
\def\sn{\textrm{sn}}
\def\ct{\textrm{ct}}
\def\cs{\textrm{cs}}
\def\re{\textrm{Re}}
\def\im{\textrm{Im}}
\def\dd{\textrm{d}}




\title{Closed Weingarten  Hypersurfaces in Warped Product Manifolds}

\author{Francisco J. Andrade {$^\dagger$}\thanks{Supported
by a PICDT scholarship.}\and Jo\~ao L. M. Barbosa
$^\ddag$\thanks{Partially supported by CNPq and PRONEX-FUNCAP.} \and
Jorge H. S. de Lira\thanks{Partially supported by CNPq and
FUNCAP-PPP.}}

\date{ }

\maketitle

\abstract{Given a compact Riemannian manifold $M$, we consider a
warped product  $\bar M = I \times_h M$  where $I$ is an open
interval in $\Rr$. We suppose that the mean curvature of the
fibers do not change sign. 
Given a positive differentiable function $\psi$ in $\bar M$, we find
a closed hypersurface $\Sigma$  which is solution of an equation of
the form $F(B)=\psi$, where  $B$ is the second fundamental form of
$\Sigma$ and $F$ is a function satisfying certain structural
properties. As examples, we may exhibit examples of hypersurfaces
with prescribed higher order mean curvature.}


\section{Introduction} \label{introduction}



Let $M^n$ be a compact Riemannian manifold and let  $I$ be an open
interval in $\mathbb{R}$. Given a positive differentiable function
$h:I\to \Rr$ we then consider the product  manifold $\bar M = I
\times M$ endowed with a warped metric
\begin{equation} \label{defmetric} \dd s^2 = \dd t^2 + h^2(t)\,\dd
\sigma^2,
\end{equation}
where   $\dd\sigma^2$ stands for the metric of $M$. We denote the
warped metric simply by $\ae \cdot , \cdot \ad$.

Given a differentiable function $z:M\to I$ its graph is defined as
the hypersurface $\Sigma$ whose points are of the form
$X(u)=(z(u),u)$ with $u\in M$. This graph is diffeomorphic with $M$
and may be globally oriented by an unit normal vector field $N$ for
which it holds that $\ae N,\partial_t\ad <0$. With respect to this
orientation, let $\lambda=(\lambda_1, \ldots, \lambda_n)$  be the
vector whose components $\lambda_i$ are the principal curvatures of
$\Sigma$, that is, the eigenvalues of the second fundamental form
$B=-\ae \dd N, \dd X\ad$ in $\Sigma$.

Let $\Gamma$ be an open convex cone with vertex at the origin in
$\Rr^n$ and containing the positive cone. Suppose that $\Gamma$ is
symmetric with respect to interchanging coordinates of its points.
Let $f$ be a positive differentiable concave function defined in
$\Gamma$. In what follows,  $f$ is supposed to be  symmetric in
$\lambda_i$ and it is required that its derivatives  satisfy $f_i>0$
in $\Gamma$.

We may define a function $F$ in the space of symmetric $n\times n$
matrices $\mathcal{S}$ setting  $F(B)=f(\lambda)$ so that  it makes
sense to write down
\[
F(B(z(u))) = f(\lambda(X(u)))
\]
when the function $z$ is supposed to be {\it admissible}, which
means that $\lambda(z(u))\in \Gamma$ for all $u\in M$. Finally,
given a positive differentiable function $\psi:\bar M\to \Rr$, it is
geometrically relevant to pose the problem of finding an admissible
function $z$ which solves the following equation
\begin{equation}
\label{eqtn} F(B(z(u)))=\psi(z(u),u), \quad u\in M.
\end{equation}
Since the second fundamental form $B$  may be written in terms of
$z$ and its first and second derivatives it happens that in
analytical terms this problem is equivalent to prove the existence
of solutions for a rather complicated fully nonlinear second order
elliptic equation. Naturally, we must impose some additional
conditions on the ambient geometry and on the structure of $f$ and
$\psi$ in order to provide a solution to (\ref{eqtn}).

Concerning the ambient geometry, we must suppose that the leaves
$M_t =\{(t,u):u\in M\}$ are mean convex with respect to the inward
unit normal vector field $-\partial_t$. This amounts to be
equivalent to the condition that
\begin{equation}
\label{hypothesis} \kappa(t)>0, \quad t\in I,
\end{equation}
where $\kappa=h'/h$. Let $\delta$ be  a strictly increasing and
continuous  function  satisfying $\delta(f)>0$ whenever $f\ge c_0$
for some positive constant $c_0$. We suppose that
\begin{equation}
\label{delta} \sum_i f_i \ge \delta(f),\quad  \sum_i f_i \lambda_i
\ge \delta(f)
\end{equation}
in points of the set
\[
\Gamma_{\mu_1,\mu_2}=\{\lambda\in\Gamma: \mu_1\ge f(\lambda)\ge
\mu_2\},
\]
where $\mu_1$ and $\mu_2$ are constants with $\mu_2\ge \mu_1>0$.
Denoting $\psi_0=\inf\psi$ we also require that
\begin{equation}
\label{limsup} \limsup_{\lambda\to \partial \Gamma} f(\lambda)\le
\bar\psi_0,
\end{equation}
for some constant $\bar \psi_0<\psi_0$.  Finally we denote
$k=f(\kappa)$. Following this notation, we state our main result.

\vspace{0.3cm}

\noindent
\begin{theorem}
\label{main}
Let $\bar M^{n+1}= I \times M^n$ be endowed with the warped metric
given by {\em (\ref{defmetric})}. Given $t_-, t_+$ with $t_-<t_+$,
consider the region $\bar M_{t_-,t_+} = \{(t,p):\; t_- \le t \le t_+
\}.$ Suppose that $f$ and $h$ satisfy the conditions  {\em
(\ref{hypothesis})-(\ref{limsup})} and suppose that $\psi$ satisfies
\begin{enumerate}

\item[a)] $\psi (t,p) > k(t)$ for $t\le t_-$,

\item[b)] $\psi (t,p) < k(t)$ for $t\ge t_+$,

\item[c)] $\partial_t \big(h\psi\big) \le 0$ for $t_-<t<t_+$.
\end{enumerate}
Then there exists  a differentiable function $z:M^n\rightarrow I$
for which
\begin{equation} \label{eq:main1}
F(B(z(u))) - \psi(z(u),u) = 0
\end{equation}
whose graph $\Sigma$ is contained in the interior of $\bar
M_{t_-,t_+}$.
\end{theorem}

\noindent Important particular cases of this theorem concern
prescribing the $r$-th mean curvatures ${n\choose r } \,H_r(\lambda)
= S_r(\lambda)$, where $S_r$ are the elementary symmetric functions
of the principal curvatures which appear in the expansion of the
characteristic polynomial of $B$. It may be seen for instance in
\cite{Trudinger} and  \cite{Spruck} that these functions fit in our
hypothesis if we consider the suitable G{\aa}rding cone.  In this
sense, the theorem  above may be viewed as an extension of
existence results found in previous contributions to the subject,
notably the works \cite{Bakelman},
\cite{TW}, 
\cite{Oliker},  \cite{Delanoe},
\cite{CNSIV}, \cite{GE}, \cite{OYYL} and \cite{QYYL}.
In these articles, it is assumed that the variation rate of $\psi$
is controlled in a certain way by the curvature of  ambient geodesic
spheres. For instance, this hypothesis in \cite{CNSIV} is stated in
terms of our notation as $\partial_t (t\psi)\le 0$ in $\bar
M_{t_-,t_+}$. Here, this hypothesis corresponds to item $(c)$ in the
statement of the theorem.

We  intend in this paper to show that the powerful elliptic tools
presented in the references above are flexible enough to be used in
a very
general geometrical setting. 
Warped products constitute a large family of Riemannian manifolds
that includes geodesic discs in space forms for suitable choices of
$I$ and $h$. Its importance as examples is pervasive in Riemannian
Geometry.

The paper is organized as follows. In Section 2, we fix notation and
present some geometric and analytic preliminaries, including the
detailed description of the problem. In Section 3 we show that under
the hypothesis of the theorem, the solutions to the problem remain
in the region $\bar M_{t_-,t_+}$. In the next section we compute
gradient and Hessian of functions which resemble the classical
height and support functions. Gradient estimates are obtained in
Section 5. The Hessian estimates proved in Section 6 are largely
inspired by the technique in \cite{QYYL}. The degree theoretical
approach to solving the problem is presented in the last section and
it is based on \cite{YYL}, \cite{OYYL} and \cite{QYYL}.




\section{Preliminaries} \label{prelim}

In the sequel, we use Latin lower case letters $i,j,\ldots$ to refer
to indices running from $1$ to $n$ and $a,b,\ldots $ to indices from
$0$ to $n$. The  Einstein summation convention is used throughout
the paper. Exceptions to these conventions will be explicitly
mentioned.

We denote the metric  (\ref{defmetric}) in $\bar M$ by $\langle
\cdot,\cdot \rangle$. The corresponding Riemannian connection in
$\bar M$ will be denoted by $\bar\nabla$. The usual connection in
$M$ will be denoted $\nabla'$. The curvature tensors in $M$ and
$\bar M$ will be denoted by $R$ and $\bar R$, respectively.

 Let $e_1, \ldots, e_n$ be an orthonormal
frame field in  $M$ and  let $\theta^1, \ldots, \theta^n$ be the
associated dual frame. The connection forms $\theta^i_{j}$ and
curvature forms $\Theta^i_{j}$ in $M$ satisfy the structural
equations
\begin{eqnarray}
& & \dd\theta^i +  \theta^i_{j}\wedge \theta^j=0,\quad \theta^i_j=-\theta^j_i, \\
& &  \dd\theta^i_{j} + \theta^i_{k}\wedge \theta^k_{j}=\Theta^i_{j}.
\end{eqnarray}
An orthonormal frame in $\bar M$ may be defined by $\bar{e}_i =
(1/h)e_i,\, 1\le i\le n,$  and $\bar{e}_{0} =
\partial/\partial t$. The associated dual frame is then
$\bar{\theta}^i = h\theta^i$ for $1\le i\le n$ and
$\bar{\theta}^{0}=\dd t$. A simple computation permits to obtain the
connection forms $\bar{\theta}^a_{b}$ and the curvature forms
$\bar{\Theta}^a_{b}$  that are given by
\begin{eqnarray}
\bar{\theta}^i_{j} &=& \theta^i_{j}, \label{omegaij}\\
\bar{\theta}_{0}^{i} &=&  (h'/h)\bar{\theta}^i, \label{omegain+1}\\
\bar{\Theta}^i_{j} &=& \Theta^i_j -(h'^2/h^2)\,\bar{\theta}^i \wedge
\bar{\theta}_j  \label{Omegaij}, \\
\bar{\Theta}_{0}^{i}
&=&({h''}/{h})\,\bar{\theta}_{0}\wedge\bar{\theta}^i,
\label{Omegain+1}
\end{eqnarray}
where  $'$ denotes the derivative with respect to $t$. Our
convention here is that
\[
\bar\theta^j_i = \ae \bar\nabla e_i,e_j\ad,\quad \bar\Theta^i_j=\ae
\bar R(\,\cdot\,,\,\cdot\,)e_j,e_i\ad.
\]
with
\[
\bar R(v,w)=\bar\nabla_v \bar\nabla_w-\bar\nabla_w\bar\nabla_v
-\bar\nabla_{[v,w]}.
\]

The frame $\bar e_a$ we just defined is adapted to the level
hypersurfaces $M_t = \{(t,p);\; p\in M \}$.  It follows from
(\ref{omegain+1}) that each fiber $M_t$ is umbilical with principal
curvatures
\begin{equation} \label{kappat}
\kappa(t) = h'(t)/h(t)
\end{equation}
calculated with respect to the {\it inward} unit normal $-\bar
e_0=-\partial/\partial t$. Notice that according our convention the
Weingarten operator for the leaves with respect to this orientation
is defined as
\[
\ae \bar\nabla e_0, e_i\ad = \bar\theta^i_0.
\]

Now, consider a smooth function $z:M\rightarrow I$. Its graph is the
regular hypersurface
\[
\Sigma = \{X(u)=(z(u),u):\;u\in M\},
\]
whose tangent space is spanned at each point by the vectors
\begin{equation} \label{tangvect1}
X_i = h\,\bar{e}_i + z_i\,\bar{e}_{0},
\end{equation}
where $z_i$ are the components of the differential $\dd z=
z_i\theta^i$. The unit vector field
\begin{equation}\label{normalvect1}
N = \frac{1}{W}\big( \sum_{i=1}^n z^i\bar{e}_i - h\bar{e}_{0}\big)
\end{equation}
is normal to $\Sigma$, where
\begin{equation} \label{wz}
W = \sqrt{h^2 + |\nabla' z|^2}.
\end{equation}
Here, $|\nabla' z|^2=z^iz_i$ is the squared norm of $\nabla
'z=z^ie_i$. The induced metric in  $\Sigma$ has components
\begin{equation}
\label{gij1} g_{ij} = \langle X_i, X_j\rangle = h^2\delta_{ij} +
z_iz_j
\end{equation}
and its inverse has components given by
\begin{equation}
\label{gij2} g^{ij}=\frac{1}{h^2}\delta^{ij}-\frac{1}{h^2W^2}z^iz^j.
\end{equation}
One easily verifies that the  second fundamental form $B$ of
$\Sigma$ with components $(a_{ij})$   is determined by
\begin{eqnarray*}
\label{bij1} a_{ij} =\langle \bar\nabla_{X_j}X_i,N\rangle=
\frac{1}{W}\big(-hz_{ij} + 2h'z_iz_j + h^2h'\delta_{ij}\big)
\end{eqnarray*}
where $z_{ij}$ are the components of the Hessian $\nabla'^2z
=\nabla'\dd z$ of $z$ in $M$.
Now, we must compute  the components  $a^i_{j}=\sum_k g^{ik}a_{kj}$
of the Weingarten map $A^\Sigma$. To simplify computations, in a
fixed point $\bar u \in M$ where $\nabla' z \neq 0$, we choose
$e_1|_{\bar u} = \nabla' z / |\nabla' z|$. We call this frame a {\it
special frame} at $\bar u$. For this choice, we obtain $\dd z =
z_1\theta^1$ at $\bar u$. Since the matrices $g_{ij}|_{\bar u}$ and
$g^{ij}|_{\bar u}$ are diagonal in a special frame, one obtains at
$\bar u$
\begin{eqnarray} \label{aij1} \nonumber
a^1_{1} &=& \frac{1}{W^3}\big(-hz_{11} + 2h'z_1^2 + h^2h'\big), \\
a^1_{i} &=& -\frac{h}{W^3}z_{1i} \quad\quad\quad\quad\quad\quad\quad\quad\quad \mbox{for}\quad 2\le i \le n, \\
\nonumber a^i_{j} &=& \frac{1}{h^2W}\big(-hz_{ij} +
h^2h'\delta_{ij}\big) \quad\quad\,\, \mbox{for} \quad  2\le i,j \le
n.
\end{eqnarray}
Special frames are quite useful for computing second and third order
covariant derivatives of $z$. By definition the Hessian of $z$ is
\begin{eqnarray}\label{nabladz}
z_{ik}\theta^k = \nabla'^2 z(e_i;\cdot)=\dd z_i -\theta^k_i z_k.
\end{eqnarray}
The third derivative of $z$ is defined by
\begin{equation}
z_{ijk}\theta^k =\nabla'^3(e_i,e_j;\cdot) =\dd z_{ij}-\theta_i^k
z_{kj}-\theta_j^k z_{ik}.
\end{equation}
Exterior differentiation of both sides in (\ref{nabladz}) gives a
Ricci identity
\begin{equation}
z_{ijk}\theta^j\wedge \theta^k=\Theta^r_i z_r
\end{equation}
and in particular (for a special frame)
\begin{equation}
\label{Ricci-1} z_{1ii}-z_{ii1}=z_{i1i}-z_{ii1}=K_i z_1,
\end{equation}
where
\begin{equation}
\label{K} K_i =\ae R(e_1,e_i)e_i,e_1\ad.
\end{equation}

Now, we consider an adapted frame field $E_0=N, E_1,\ldots, E_n$ in
some open set in $\Sigma$. Representing by $\omega^a$ its dual
forms, by $\omega^a_{b}$ its connection forms and by
$\bar\Omega^a_{b}$ its curvature forms, we have the following
relations:
\begin{eqnarray} \label{eq604}
& & d\omega^i +  \omega^i_{j}\wedge \omega^j=0,\quad \omega^i_j=-\omega^j_i, \\
& & d\omega^i_{j} + \omega^i_{k}\wedge\omega^k_{j} = \Omega^i_{j},
\end{eqnarray}
where $\Omega^i_j$ are  the curvature forms for $\Sigma$. Since
$\Sigma$ is a hypersurface of $\bar{M}$ then the Gauss equation
reads off as
\begin{equation} \label{eq605}
\Omega^i_{j} = \bar\Omega^i_{j} - \omega^i_{0}\wedge\omega^{0}_j
\end{equation}
The coefficients $a_{ij}$ of the second fundamental form are given
by Weingarten equation
\begin{equation}
\label{Weingarten} \omega^{0}_i = a_{ij}\,\omega^j.
\end{equation}
In the sequel, one
indicates the covariant derivative in $\Sigma$  by $\nabla$ and by a
semi-colon. Remember that
\begin{eqnarray}  \label{eqn01}
\nabla a_{ij} &=& \dd a_{ij} - a_{kj}\,
\omega^k_{i} -  a_{ik}\,\omega^k_{j} =  a_{ij;k}\,\omega^k \\
\nonumber \nabla a_{ij;k} &=& \dd a_{ij;k} - a_{mj;k}\,\omega^m_{i}
-
a_{im;k}\,\omega^m_{j} -  a_{ij;m}\,\omega^m_{k} \\
\label{eqcl02} &=&  a_{ij;km}\,\omega^m
\end{eqnarray}
The Codazzi equation is a commutation formula for the first
derivative of $a_{ij}$ and it is obtained by differentiating
(\ref{Weingarten}):
\begin{eqnarray}
\label{Codazzi} a_{ij;k}\omega^j\wedge \omega^k=\bar\Omega_{i}^{0}.
\end{eqnarray}

We also prove using the preceding notation a very useful Ricci
identity.

\begin{lemma} \label{prop64}
Let $\bar{X}$ be a point of $\Sigma$ and $E_0 =N, E_1,\ldots, E_n$
be an adapted frame field such that each $E_i$ is a principal
direction and $\omega^k_i=0$ at $\bar X$. Let $(a_{ij})$ be the
second quadratic form of $\Sigma$. Then, at the point $\bar{X}$, we
have
\[
a_{ii;11}-a_{11;ii}=a_{11}a_{ii}^2-a_{11}^2a_{ii}+\bar
R_{i0i0}\,a_{11} -\bar R_{1010}\,a_{ii} +\bar R_{i1i0;1}-\bar
R_{1i10;i}.
\]
\end{lemma}



The frame field $E_a$ may be obtained from the adapted frame field
$N,X_1,\ldots, X_n$ by Gram-Schmidt procedure. Since this last frame
depends only on $z$ and $\nabla' z$, we may conclude that components
of $\bar R$ and $\bar\nabla \bar R$ calculated in terms of the frame
$E_a$ depend only on $z$ and $\nabla' z$.

\subsection{The prescribed curvature equation}

Now we formulate the existence problem analytically. We consider $f$
and $\Gamma$ as defined in Section 1. Then, given the second
fundamental form $(a_{ij})$ in $\Sigma$ we define
\[
F\big((a_{ij})\big)=f(\lambda_1,\ldots,\lambda_n),
\]
where $\lambda_i$ are the eigenvalues of $(a_{ij})$ calculated with
respect to the induced metric $(g_{ij})$. It is convenient to denote
the vector of principal curvatures $(\lambda_1,\ldots, \lambda_n)$
by $\lambda$. Admissible functions are those ones for which
$\lambda$ always lies in $\Gamma$. We may consider $F$ as a map from
$\mathcal{S}\times\Rr^n\times\Rr$ into $\mathbb{R}$  in the
variables $z_{ij}$, $z_i$ and $z$.

Thus our problem is to find $\Sigma$, graph of an admissible
function, so that
\[
F\big(a_{ij}(z(u))\big)=\psi(z(u),u),\quad u\in M,
\]
for some prescribed positive function $\psi$.  We recall that is
required  that $f$ satifies
\begin{equation}
\label{f-elliptic} f_{i}=\frac{\partial f}{\partial \lambda_i}>0
\end{equation}
and that $f$ is concave what implies that
\begin{equation}
\label{f-concavity} \sum_i f_i \lambda_i \le f.
\end{equation}
We also assume the condition (\ref{delta}) and then we prove using
the assumption (\ref{limsup}) and following \cite{CNSIV} that
\[
\sum_i \lambda_i\ge \delta
\]
for $\lambda\in \Gamma$ such that $f(\lambda)\ge \psi_0$. In fact,
the set
\[
\Gamma_{\psi}=\{\lambda\in\Gamma : f(\lambda)\ge \psi_0\}
\]
is closed in $\Rr^n$, convex and symmetric. Thus the closest point
in $\Gamma_\psi$ to the origin is of the form $(\lambda_0,\ldots,
\lambda_0)$. This geometric fact implies that  any $\lambda\in
\Gamma_\psi$ is located above the hyperplane
\begin{equation}
H=\Big\{\lambda\in \Rr^n: \sum_i \lambda_i =n\lambda_{0}\Big\}.
\end{equation}
Hence, any $\lambda\in \Gamma_\psi$ is necessarily contained in the
convex part of the cone $\Gamma$ which is above $H$. This implies
that upper estimates for $\lambda$ imply automatically lower
estimates.

We proceed by stating some useful analytical properties of $F$.
Notice that $F$ is differentiable whenever $f$ is.   We denote first
and second derivatives of $F$ respectively by
\begin{equation}
F^{ij} =  \frac{\partial F}{\partial a_{ij}}\quad \textrm{  and
}\quad F^{ij,kl}=\frac{\partial^2 F}{\partial a_{ij}\partial
a_{kl}}.
\end{equation}
These derivatives may be easily calculated if we assume that the
matrix $(a_{ij})$ is diagonal with respect to $(g_{ij})$, due to the
following lemma.

\begin{lemma}
\label{deriveSm} If $(a_{ij})$ is diagonal at $\bar X$ then the
matrix $(F^{ij})$ is also diagonal with positive eigenvalues $f_i$.
Moreover, $F$ is concave and its second derivatives are given by
\begin{equation}
F^{ij,kl}\eta_{ij}\eta_{kl}=\sum_{k,l}f_{kl}\eta_{kk}\eta_{ll}+\sum_{k\neq
l} \frac{f_k-f_l}{\lambda_k-\lambda_l}\eta_{kl}^2.
\end{equation}
Finally one has
\begin{equation}
\frac{f_i-f_j}{\lambda_i-\lambda_j}\le 0.
\end{equation}
These expressions must be interpreted as limits in the case of
multiple eigenvalues of $(a_{ij})$.
\end{lemma}

\noindent The terms $F^{ij}$ are components of a rank two
contravariant tensor. Thus one has
\[
F^{ij}a_{ij}= F^i_j a^j_i
\]
and if the matrix $(g_{ij})$ is assumed to be diagonal at $\bar X$,
then $(F^i_j)$ is also diagonal at that point.

\section{Height estimates}

Now,  we consider, for each $s$, $0\le s \le 1$, the map
\begin{equation} \label{defPsi}
\Psi(s,t,u) = s\psi(t,u) + (1-s)\phi(t)k(t),
\end{equation}
where $k(t)=f(\kappa(t))$ and $\phi$ is a positive real function
defined in $I$, which satisfies the following conditions:
\begin{enumerate}
\item[a)] $\phi>0$,

\item[b)] $\phi(t)>1 \quad \mbox{for} \quad t\le t_-$,

\item[c)]$\phi(t)<1 \quad \mbox{for} \quad t\ge t_+$,

\item[d)] $\phi'(t)<0$.
\end{enumerate}
These conditions imply the existence of a unique point $t_0 \in
(t_-, t_+)$ such that $\phi(t_0) = 1$. Combining the conditions
above on $\phi$ and the hypothesis $(a)$ and $(b)$ in the statement
of the theorem, one proves

\vspace{3mm}

\begin{lemma} \label{aboutfamily}
For $\psi$ as in Theorem 1, $\phi$ as prescribed above and the
function $\Psi$ defined in {\em (\ref{defPsi})}, the following
statements are true:
\begin{enumerate}
\item[i)] $\Psi(1,t,u) = \psi(t,u) \quad \mbox{and} \quad
\Psi(0,t,u)=\phi(t)k(t),
$

\item[ii)] $\Psi(s,t,u)>0$,

\item[iii)] $\Psi(s,t,u)> k(t) \quad \mbox{for}\quad t\le t_-,$

\item[iv)] $\Psi(s,t,u)< k(t)
\quad \mbox{for} \quad t\ge t_+.$
\end{enumerate}
\noindent Furthermore, it is always possible to choose $\phi$
satisfying the prescribed conditions such that:
\begin{enumerate}
\item[v)]  $\frac{\partial}{\partial t}\Psi(s,t,u) + \kappa(t)
\Psi(s,t,u)<0.$
\end{enumerate}
\end{lemma}

\vspace{3mm}

\noindent For $0\le s\le 1$, consider  the family of equations
\begin{equation} \label{eqfamily}
\Upsilon(s,z)= F(a_{ij}(z)) - \Psi(s,z,u) = 0, \quad z = z(u).
\end{equation}
Notice that the constant function $t=t_0$ is solution to the problem
corresponding to $s=0$. We denote it by $z_0$.

We are able to prove  $C^0$ bounds uniform with respect to the
parameter $s$ of this homotopy. More precisely, one proves

\begin{prop} \label{prop-C0}
Suppose that   $\psi$ satisfies  the conditions {\em (a)} and {\em
(b)} in Theorem 1.
If $z\in C^2(M)$ is a solution of the equation $\Upsilon(s,z)=0$ for
a given $0\le s\le 1$, then  
\begin{equation}
t_-<z(u)<t_+, \quad u\in M.
\end{equation}
\end{prop}

\vspace{0.3cm}

\noindent{\it Proof:}  Let $\bar u$ be a point of maximum for the
function $z(u)$. This exists by the compactness of $M$. Let's assume
that $z(\bar u) \ge t_+$. Consider then the leaf $M_{z(\bar u)}$ and
represent by $\Sigma$ the graph of $z$. Observe that $\Sigma$ and
$M_{z(\bar u)}$ are tangent at $(z(\bar u),\bar u)$. Furthermore,
with respect to the inwards normal vector common to both
hypersurfaces at this point, $\Sigma$ lies above $M_{z(\bar u)}$.
But then the principal curvatures of $\Sigma$ at this point are
greater than or equal to $\kappa(z(\bar u))$. Thus by the fact that
$f$ has positive derivatives one concludes that
\[
F(a_{ij}(z))\ge k(z (\bar u))
\]
at $(z(\bar u),\bar u)$ what is in contradiction with $(iv)$ of
Lemma (\ref{aboutfamily}). Hence $z(\bar u) < t_+$. Working in a
similar way with the minimum $\hat u$ of $z(u)$ one concludes that
$z(\hat u)> t_-$.


\vspace{0.3cm}

Now, we prove the following uniqueness result.

\begin{prop}
Fixed $s=0$ there exists an unique admissible solution $z_{0}$ of
the equation $\Upsilon(0,z)=0$, namely $z_0 =t_0$ where $t_0$
satisfies $\phi(t_0)=1$.
\end{prop}

\vspace{0.3cm} \noindent{\it Proof.} That $z_0$ is solution to this
problem follows from
\[
\Upsilon(0,z_{0})= F(a_{ij}(z_{0})))-k(t_0)=f(\kappa(t_0))-k(t_0)=0.
\]
Let  $\bar{z}$ be an admissible solution of $\Upsilon(0,z)=0$. This
means that
\[
F(a_{ij}(\bar{z})) -\phi(\bar{z})k(\bar{z})=0.
\]
Now, let $\bar{u}\in M$ a minimum point of $\bar{z}$. At this point,
one has  $\nabla'\bar{z}=0$ and $\nabla'^{2} \bar z$ is
positive-definite. Since $\kappa=\frac{h'}{h}$ one computes
explicitly at $\bar u$
\[
a^i_{j}=g^{ik}a_{kj}=-\frac{1}{h^2}\sigma^{ik}\bar
z_{kj}+\frac{h'}{h}\delta^i_j
\]
Therefore if we consider a local frame around $\bar u$ which is
orthonormal at $\bar u$ and which diagonalizes $\nabla'^{2}\bar z$
at this point one obtains
\[
a^i_j(\bar z(\bar u))\le \kappa(\bar z(\bar u))\delta^i_j
\]
and since $f$ is increasing with respect to its arguments
\begin{eqnarray*}
\phi(\bar{z}(\bar{u}))k(\bar{z}(\bar{u}))=
F(a_{ij}(\bar{z}(\bar{u})))\le f(\kappa(\bar z(\bar u))) = k(\bar
z(\bar u))=\phi(t_0)k(\bar z(\bar u)).
\end{eqnarray*}
Hence, since  $\phi$  is a decreasing function one concludes from
the choice of $\bar u$ as a minimum point that
\[
\bar{z}(u)\geq\bar{z}(\bar{u})\geq t_{0},
\]
for all $u\in M$. In a similar way, one proves that
\[
\bar{z}(u)\leq t_{0}
\]
for all $u\in M$. Thus, one gets $z=z_0$. This finishes the proof.

\section{Height and support functions}

As before, let $\Sigma$ be the graph of $z$. We start by considering
the functions $\tau:\Sigma\to \mathbb{R}$ and $\eta:\Sigma\to
\mathbb{R}$ given by
\begin{eqnarray}\label{t}
\tau
=-h \ae N, \bar e_0\ad
\quad \textrm{and}\quad \eta=- \int h\, \dd t.
\end{eqnarray}
The following formulae will be useful later.
\begin{lemma} The gradient vector fields of the functions
$\tau$ and $\eta$ are
\begin{eqnarray}\label{dn}
\nabla\eta  &=& -  h\,\bar e_0^T,\\
\label{dt}\nabla\tau &=&- A^\Sigma(\nabla \eta),
\end{eqnarray}
and its Hessian forms calculated with respect to given vector fields
$V,W$ in $\Sigma$ are
\begin{eqnarray}
\label{ddn} \nabla^2\eta(V,W) &=& \tau B(V,W)-h'\ae V,W\ad,\\
 \label{ddt}
\nabla^2\tau(V, W) & = &- \ae\nabla_{\nabla\eta}A^\Sigma
V,W\ad-\ae\bar R(\nabla\eta,W)V,N\ad \nonumber\\
& -&\tau \ae A^\Sigma V, A^{\Sigma}W\ad+ h'\ae A^\Sigma V, W\ad,
\end{eqnarray}
Here, $\bar e_0^T$ denotes the tangential projection of the vector
field $\bar e_0$.
\end{lemma}


\noindent {\it Proof.} To simplify the calculations, we consider a
local  orthonormal frame $e_a$  around a point $\bar u$ of $M$ and
the associated adapted frame field $N,E_1,\ldots,E_n$ along $\Sigma$
so that $\bar\nabla E_i|_{X(\bar u)}=0$. Using (\ref{tangvect1}) one
has
\begin{eqnarray*} \dd \eta = - h\,\dd z = -h\ae \dd X,\bar
e_0\ad = -h \ae \bar e_0^T, \omega^i E_i\ad
\end{eqnarray*}
and
\begin{eqnarray*}
\dd \tau &=& -\dd h \ae N, \bar e_0\ad -h \ae\bar\nabla N,\bar
e_0\ad -h\ae
N,\bar\nabla \bar e_0\ad\\
&  = &-h'\bar\theta^0 \ae N, \bar e_0\ad +h \ae a^j_i
E_j\omega^i,\bar e_0\ad -h\ae N,\bar\theta^i_0 \bar e_i\ad\\
&  = & h \ae a^j_i E_j\omega^i,\bar e_0\ad-h'\bar\theta^0 \ae N,
\bar
e_0\ad  -h'\ae N,\bar\theta^i \bar e_i\ad\\
&  = & h \ae A^\Sigma(E_i),\bar e_0^T\ad\omega^i-h' \ae
N,\bar\theta^0 \bar e_0+\bar\theta^i \bar e_i\ad.
\end{eqnarray*}
Thus since $A^\Sigma$ is self-adjoint and  $\dd X = \bar\theta^0\bar
e_0 +\bar\theta^i\bar e_i$, one gets
\begin{equation}
\dd \tau =h \ae A^\Sigma(\bar e_0^T), \omega^i E_i\ad.
\end{equation}
Therefore we conclude that
\begin{equation}
\nabla \eta = -h \, \bar e_0^T, \quad \nabla \tau =-A^\Sigma(\nabla
\eta).
\end{equation}
Since $\tau_i = h \ae a_i^j E_j, \bar e_0^T\ad$, one computes, using
$\nabla_{E_k}E_j|_{X(\bar u)}=0$,
\begin{eqnarray*}
 \tau_{i;k} &=& h_k \ae a_i^j E_j, \bar e_0\ad  +h \ae a_{i;k}^j
E_j, \bar e_0\ad + h \ae a_i^j \bar\nabla_{E_k}E_j, \bar e_0\ad
\\
& + &  h \ae a_i^j E_j, \bar\nabla_{E_k}\bar e_0\ad \\
& = & h' \ae a_i^j E_j, \bar\theta^0(E_k)\bar e_0\ad+h \ae a_{i;k}^j
E_j, \bar e_0\ad + h  a_i^j a_{kj} \ae N, \bar e_0\ad \\
&  + & h' \ae a_i^j
E_j, \bar\theta^i(E_k)\bar e_i\ad\\
&  = & h' \ae a_i^j E_j, E_k\ad+h \ae a_{i;k}^j E_j, \bar e_0\ad + h
a_i^j a_{kj} \ae N, \bar e_0\ad\\
&  = &  h' a_{ik}-a_{i;k}^j\eta_j -\tau a^j_i a_{kj}.
\end{eqnarray*}
where we used again that $\dd X = \bar\theta^0 \bar e_0
+\bar\theta^i \bar e_i$ and that $\eta_k = - h \ae \bar e_0,
E_k\ad$. Hence, one gets from Codazzi's equation
\begin{eqnarray*}
\nabla^2 \tau(V,W) & = & 
h'\ae A^{\Sigma}V,W\ad-\ae (\nabla_W A^{\Sigma})\nabla \eta, V
\ad -\tau \ae A^\Sigma V, A^\Sigma W\ad\\
& = & h'\ae A^{\Sigma}V,W\ad-\ae (\nabla_{\nabla \eta} A^{\Sigma})
W, V \ad -\ae \bar R(\nabla\eta,W)V,N\ad\\
&- &\tau \ae A^\Sigma V, A^\Sigma W\ad.
\end{eqnarray*}
Finally, it follows from the expression $\eta_i = - h\ae E_i, \bar
e_0\ad$ that
\begin{eqnarray*}
\eta_{i;k} & = & -h_k \ae E_i, \bar e_0\ad - h\ae
\bar\nabla_{E_k}E_i, \bar e_0\ad - h \ae E_i, \bar\nabla_{E_k}\bar
e_0\ad \\
& = &-h'  \ae E_i, \bar\theta^0(E_k)\bar e_0+\bar\theta^i(E_k)\bar
e_i\ad - h\,a_{ik} \ae N, \bar e_0\ad \\
& =  & -h' g_{ik} +\tau a_{ik}.
\end{eqnarray*}
Thus we obtain
\begin{eqnarray}
\nabla^2 \eta(V,W) = -h' \ae V,W\ad + \tau \ae A^\Sigma V, W\ad.
\end{eqnarray}
This finishes proving the lemma. \vspace{0.3cm}

One estimates the derivatives of $\eta$ and $\psi$ as follows. In
the sequel $\nabla_i$ and $\nabla_{ij}$ denote covariant derivative
in $\Sigma$ calculated with respect to a frame adapted to $\Sigma$.

\begin{lemma} The functions $\eta$ and $\psi$ satisfy the following
estimates
\begin{eqnarray}
|\nabla\eta|\le C,\quad
|\nabla\psi|\le C,\quad 
|\nabla^2\psi|\le C 
\label{des3}
\end{eqnarray}
where $C$ are constants depending on $\psi$, $\nabla'\psi$,
$\nabla'^2\psi$ and on $C^0$ and $C^1$ bounds for $z$.
\end{lemma}

\noindent {\it Proof.} The first estimate follows from the $C^0$ and
$C^1$ estimates for $z$. In fact, one has $\eta_i = -hz_i$.
In order to prove the remaining estimates, we observe that
$$\nabla_{i}\psi=X_i(\psi)=e_i(\psi)+z_i e_0(\psi)=:\psi_{i}+z_i\psi_{z}.$$
Thus, using (\ref{gij2}) and denoting $\psi_i =e_i(\psi)$ and so on
we have
\begin{eqnarray*}
|\nabla \psi|^2 & = & g^{ij}X_i(\psi)X_j(\psi)= g^{ij}\big(\psi_i +
z_i\psi_z\big)\big(\psi_j+z_j\psi_z\big)\\
& = &
\frac{1}{h^2}\Big(\delta^{ij}\psi_i\psi_j-\frac{z^iz^j}{W^2}\psi_i\psi_j
+\delta^{ij}z_i\psi_j\psi_z
-\frac{z^iz^j}{W^2}z_i\psi_j\psi_z\\
& + &\delta^{ij}z_j\psi_i\psi_z
-\frac{z^iz^j}{W^2}z_j\psi_i\psi_z+\delta^{ij}z_iz_j\psi_z^2
-\frac{z^iz^j}{W^2}z_i
z_j\psi_z^2\Big)\\
& = &
\frac{1}{h^2}\Big(|\nabla'\psi|^2-\frac{1}{W^2}\ae\nabla'\psi,\nabla'z\ad^2+2\psi_z\ae
\nabla'\psi,\nabla'z\ad\\
&- & 2\frac{\psi_z}{W^2}| \nabla'z|^2\ae \nabla'\psi,\nabla'z\ad +
\psi_z^2|\nabla'z|^2-\frac{\psi_z^2}{W^2}|\nabla'z|^4\Big)\\
&\le &  C(|z|_1,|\psi|_1,|\psi_z|).
\end{eqnarray*}
In a similar way (replacing $\psi$ by $\psi_t=\psi_z$) we prove that
\begin{equation}
|\nabla\psi_z|\le C.
\end{equation}
One has
\[
X_iX_j(\psi)=X_i\big(\psi_j+z_j\psi_z\big)=\psi_{i,j}+z_{i,j}\psi_z+z_j\psi_{zi}+z_i\psi_{zj}+z_iz_j\psi_{zz},
\]
where $\psi_{i,j}=e_ie_j(\psi)$ and $z_{i,j}=e_ie_j(z)$. We then
choose a geodesic frame $e_a$ around $\bar u\in M$. In this case it
holds that $z_{i,j}=\nabla'_{ij}z=z_{ij}$ at $\bar u$.
Now using the fact that $\bar \theta^a_b=0$ at $\bar u$, we obtain
\begin{eqnarray*}
\bar\nabla_{X_j}X_i &=& \big(\dd z_i(X_j)+h'h\theta^i(X_j)\big)\,e_0
+
\frac{h'}{h}\,\big(\delta^k_i\theta^0(X_j)+z_i\theta^k(X_j)\big)\,e_k\\
& = & \big(z_{ij}+h'h\delta_{ij}\big)\,e_0 +
\frac{h'}{h}\,\big(z_je_i+z_ie_j\big)
\end{eqnarray*}
which implies that
\begin{eqnarray*}
& &  \ae \nabla_{X_j}X_i, \nabla\psi\ad
 = \big(z_{ij}+h'h \delta_{ij}\big)\psi_z + \frac{h'}{h}\big(z_j
\psi_i+z_i\psi_j\big).
\end{eqnarray*}
Hence, one obtains
\begin{eqnarray*}
\nabla_{ij}\psi &=& \ae \nabla_{X_j}\nabla \psi, X_j\ad
=\psi_{ij}+z_j\psi_{zi}+z_i\psi_{zj}+z_iz_j\psi_{zz}\\
& - & h'h \delta_{ij}\psi_z - \frac{h'}{h}z_j
\psi_i-\frac{h'}{h}z_i\psi_j.
\end{eqnarray*}
Therefore we conclude that
\[
|\nabla^2 \psi|\le C(|z|_1,|\psi|_2, |\psi_z|_1).
\]
This finishes the proof of the lemma.


\section{Gradient estimate} \label{gradient}

In this section, we prove {\it a priori} global estimate for the
first derivatives of $z$.

\begin{prop}
\label{est-nabla} Under the hypothesis of Theorem 1,  if $z(u)$ is a
solution of equation {\em (\ref{eqfamily})} for some fixed $0\le
s\le 1$, then $|\nabla' z|< C $, where $C$ is a constant that
depends only on $t_-$, $t_+$ and $\psi$.
\end{prop}

\noindent {\it Proof.} We present the proof for $s=1$. There is no
essential changes for $0\le s<1$.

Set $\chi(z) = |\nabla' z|e^{A z}$, where $A$ is a positive constant
to be chosen later on. Let $\bar u$ be a point where $\chi$ attains
its maximum. If $\chi(\bar u) = 0$ then $|\nabla' z|\equiv 0$ and so
the result is trivial. Hence, we are going to assume that $\chi(\bar
u)>0$. Thus we may define the function $ \ln \chi(z) = \ln |\nabla'
z| + Az$ which also attains its maximum at $\bar u$. Hence, fixing a
special frame in some neighborhood of $\bar u$ one has
\begin{eqnarray*}
0 & = & \chi_i = \frac{1}{|\nabla'z|}e^{Az}z_{ik}z^k+
Ae^{Az}|\nabla'z|z_i\\
& = & e^{Az}z_{i1}+ Ae^{Az}z^1z_i,
\end{eqnarray*}
which implies by the symmetry $z_{i1}=z_{1i}$ of the Hessian form
that
\begin{equation}
\label{eq-z11} z_{11}=-Az_1^2,\quad z_{1i}=0, \quad i>1.
\end{equation}
where we used the fact that $z_i|_{\bar u}=0$ for $i\neq 1$.
Substitution of this into (\ref{aij1}) yields $a_{1i} = 0$ for
$i>1$. This implies that the direction $e_1$ at $\bar u$ is
principal. Then, we may rotate the other vectors $e_2, \ldots, e_n$
so that they are also principal at $\bar u$. With this choice we
have $a_{ij} = 0$ for $i\ne j$ ate $\bar u$. As a consequence of
this, one sees from (\ref{aij1}) that $z_{ij}(\bar u) = 0$ for $i\ne
j$. Thus, the Hessian of $z$ is diagonal at $\bar u$.

Differentiating again the function $\chi$ at $\bar u$, one obtains
(no summation over the index $i$)
\begin{eqnarray*}
0 &\ge & \chi_{i;i}=Ae^{Az}z_iz_{i1}+A^2 e^{Az}z_1 z_i^2\\
& + & e^{Az}\left(-\frac{1}{z_1}z_{i1}^2
+z_{i1i}+\frac{1}{z_1}z_{ii}^2+Az_iz_{i1}+Az_1z_{ii}\right).
\end{eqnarray*}
Hence, one concludes  from this inequality that
\begin{eqnarray}  
z_{111} +A^2 z_1^3 + 3Az_1 z_{11} \le 0, &&\\
\label{eq-B}z_{i1i} +\frac{1}{z_1}z_{ii}^2+ Az_1 z_{ii} \le 0. &&
\end{eqnarray}
Combining the first inequality just above and (\ref{eq-z11}) gives
\begin{equation}
\label{eqA1} z_{111}-2A^2 z_1^3\le 0.
\end{equation}
From (\ref{eq-B}) and (\ref{Ricci-1}) one gets
\begin{equation} \label{eqB2}
z_{ii1} \le - \frac{z_{ii}^2}{z_1} - Az_1 z_{ii} - K_iz_1 \quad
\mbox{for } i>1\,.
\end{equation}

Now we can start putting all this information together to obtain the
desired estimate. We start by taking the derivative of equation
(\ref{eqfamily}) with respect to the direction $e_1$. Using the fact
that the matrix $(a^i_{j})$ is diagonal at $u_0$ and the remarks
just after Lemma \ref{deriveSm}, we obtain:
\begin{equation} \label{principal}
\sum_{i=1}^n F^{i}_i\frac{\partial a^i_{i}}{\partial z_1}z_{11} +
\sum_{i=1}^n F^{i}_i\frac{\partial a^i_i}{\partial z}z_1 = \psi _z
z_1 - \sum_{i=1}^n F^{i}_i\frac{\partial a^i_{i}}{\partial
z_{ii}}z_{ii1}.
\end{equation}
Taking derivatives of $a^i_{i}$, using (\ref{aij1}) we obtain
\begin{eqnarray*}
\frac{\partial a^1_{1}}{\partial z_{11}} & =&  - \frac{h}{W^3},\\
\label{eqd111} \frac{\partial a^1_{1}}{\partial z_1} & = &  -
\frac{3z_1}{W^2}a^1_{1} + \frac{4z_1h'}{W^3},\\
\nonumber  \frac{\partial a^1_{1}}{\partial z}  & = & \left(\frac{h'
}{h} - \frac{3hh'}{W^2}\right)a^1_{1}
+ \frac{2}{h W^3}\big(hh''-h'^2\big)z_1^2\\
&+&\frac{1}{W^3}\big(hh'+h^2h''\big)
\end{eqnarray*}
and for $i>1$
\begin{eqnarray*}
\label{eqdiiii} \frac{\partial a^i_{i}}{\partial z_{ii}} & = &  -
\frac{1}{ h W},\\
\label{eqdii1} \frac{\partial a^i_{i}}{\partial z_1} & = &  -
\frac{z_1}{W^2}a^i_{i},\\
\label{eqdii} \frac{\partial a^i_{i}}{\partial z} & = &
-h'\left(\frac{h}{W^2} + \frac{1}{h}\right)a^i_{i} + \frac{(hh')'}{h
W}\,.
\end{eqnarray*}
Replacing this into (\ref{principal}) and using (\ref{eq-z11}) and
rearranging terms  yields
\begin{eqnarray} \nonumber  \label{eqded04}
& & z_1\left(\frac{3Az_1^2}{W^2} + \frac{h'}
{h} - \frac{3hh'}{W^2}\right)F^{1}_1a_{1}^1  \\
\nonumber & &  + z_1 \left(-\frac{4Az_1^2h'}{W^3}
+\frac{2}{hW^3}(hh''-h'^2)z_1^2 +
\frac{1}{W^3}(hh'+h^2h'')\right)F^{1}_1 \\
& & \nonumber+ z_1\left(\frac{Az_1^2}{W^2} - h'\left(\frac{h}{W^2} +
\frac{1}{h}\right)\right)\sum_{i>1}F^{i}_ia^i_{i}+ z_1\frac{(hh')'}{hW}\sum_{i>1}F^{i}_i\\
& & = \psi_z z_1 + F^{1}_1\frac{h}{W^3}z_{111} + \sum_{i>1}
F^{i}_i\frac{1}{h W}z_{ii1 }.
\end{eqnarray}
Using (\ref{eqA1}) and (\ref{eqB2}) we can estimate the right hand
side of (\ref{eqded04}) by
\begin{eqnarray}\nonumber \label{eqded01}
& & \textrm{RHS} \le   \psi_z z_1 + F^{1}_1\frac{2A^2h z_1^3}{W^3} -
Az_1\sum_{i>1}F^{i}_i\frac{z_{ii}}{ h W} - \frac{K_i z_1}{ h
W}\sum_{i>1}
F^{i}_i\\
& & \nonumber\le\psi_z z_1 + Az_1\sum_i F^i_i a^i_i +
F^{1}_1\left(\frac{A^2h z_1^3}{W^3}-
Az_1\frac{h'}{W^3}(2z_1^2 + h^2) \right) \\
& & -\left( \frac{A h'z_1}{W} + \frac{K_i z_1}{hW}\right)\sum_{i>1}
F^{i}_i,
\end{eqnarray}
where we used the expressions of $a_{1}^1$ and $a_{i}^i$ given in
(\ref{aij1}) and the fact that $F^{i}_i>0$.

Transposing the term in $\sum_{i>1}F^{i}_i$ from the right hand side
in (\ref{eqded01}) to the left hand side of the equation
(\ref{eqded04}), and adding it with the one that was already there
and finally choosing $A$ so that
\begin{equation}
\label{A} Ah'h+(h'h)'+\min_i K_i >0
\end{equation}
results, by the fact that $h'
>0$, in a positive term that can be discarded. Notice that $K_i=\ae R(e_1,e_i)e_1,e_i\ad$
does not depend on derivatives of $z$. This and the fact that $h$
and its derivatives are uniformly bounded in the annulus $\bar
M_{t_-,t_+}$ show that we may choose any $A\ge A_0$ for some $A_0$
which depends only on $t_-, t_+$ and $|z|_0$.

We may estimate the left hand side of the inequality resulting from
(\ref{eqded04}) after these manipulations  as
\begin{eqnarray} \nonumber \label{eqded06}
\textrm{LHS} &\ge &  z_1\left(\frac{Az_1^2}{W^2}  -
\frac{hh'}{W^2}-\frac{h'}{h}\right)\sum_i F^i_i a^i_i \\
&  + & z_1\left(\frac{2Az_1^2}{W^2} +
\frac{2h'}{h} - \frac{2hh'}{W^2}\right)F^{1}_1a_{1}^1 \nonumber \\
&  + & \frac{z_1}{W^3}\bigg(-4Ah'z_1^2 + \frac{2z_1^2}{h}\big(hh''
- h'^2\big)\nonumber \\
& + & hh' + h^2h''\bigg)F^{1}_1.
\end{eqnarray}
Transpose the term with $F^{1}_1$ from the right hand side in
(\ref{eqded01}) to the right hand side in (\ref{eqded06}) and add it
to the one that exists there. Transpose the term in $\sum_i F^i_i
a^i_i$ from the right hand side in (\ref{eqded06}) to the right hand
side of the inequality  (\ref{eqded01})  obtaining
\begin{eqnarray}  \label{eqded03B}
\textrm{RHS} \le  \psi_z z_1+Az_1\sum_i F^i_i a^i_i -
z_1\big(\frac{Az_1^2}{W^2}-\frac{hh'}{W^2}-\frac{h'}{h}\big)\sum_i
F^i_i a^i_i.
\end{eqnarray}
For the left hand side we obtain
\begin{eqnarray} \nonumber \label{eqded06B}
\textrm{LHS}  &\ge &   2z_1\left(\frac{Az_1^2}{W^2} + \frac{h'}{h} -
\frac{hh'}{W^2}\right)F^{1}_1a_{1}^1 + \frac{z_1}{W^3}\bigg(\frac{2z_1^2}{h}(hh''-h'^2)\\
&  + & hh' + h^2h'' +  Ah'(-2z_1^2+h^2) - A^2z_1^2h \bigg)F^{1}_1.
\end{eqnarray}
Thus, replacing in (\ref{eqded06B}) the expression for $a^1_1$ in
(\ref{aij1}) and gathering the resulting expression to
(\ref{eqded03B}), one gets
\begin{eqnarray}
\label{desfinal} & & \frac{2z_1}{W^3}\left(\frac{Az_1^2}{W^2} +
\frac{h'}{h} -
\frac{hh'}{W^2}\right)\left(Ahz_1^2+2h'z_1^2+h^2h'\right)F^{1}_1\nonumber\\
& & +
\frac{z_1}{W^3}\bigg(\frac{2z_1^2}{h}(hh''-h'^2)+ hh' + h^2h'' +  Ah'(-2z_1^2+h^2) - A^2z_1^2h \bigg)F^{1}_1\nonumber\\
& & \label{pol-1}\le  \psi_z z_1+Az_1\sum_i F^i_i a^i_i -z_1
\big(\frac{Az_1^2}{W^2}-\frac{hh'}{W^2}-\frac{h'}{h}\big)\sum_i
F^i_i a^i_i.
\end{eqnarray}
Observe that in (\ref{desfinal}) all coefficients of $F^1_1$  have
uniform lower bounds and moreover that the first term in the left
hand side of (\ref{desfinal}) is nonnegative. Thus, it is possibel
to consider this inequality as polynomial in $A$ writing it as
\begin{eqnarray}
\label{pol} F^1_1 \big(aA^2 +bA +c\big) &\le & \psi_1+\psi_z
z_1+Az_1\sum_i F^i_i a^i_i \nonumber\\
&-&z_1\sum_i F^i_i a^i_i
\big(\frac{Az_1^2}{W^2}-\frac{hh'}{W^2}-\frac{h'}{h}\big),
\end{eqnarray}
where  $a,b$ e $c$ are coefficients uniformly bounded in terms of
the functions  $h$, $h'$ e $h''$. Thus, we must consider two cases.
First, we suppose that  $F^1_1$ is uniformly bounded from zero,
i.e., that there exists a constant $C>0$ such that  $F^1_1\ge C$ em
$\Sigma$. In this case, the coefficient
\begin{equation}
a = \frac{hz_1^3}{W^5}(z_1^2 - h^2)F^{1}_1
\end{equation}
is necessarily nonpositive, since $A$ may be chosen arbitrarily
large in (\ref{pol}). Thus, it follows that  $z_1(\bar u) \le
h(z(\bar u))$ and therefore $z_1(\bar u) < h(t_+)$.

The other possibility is that $F^1_1$ has no strictly positive lower
bound. In this case, it is convenient to write the left hand side in
(\ref{desfinal}) as
\begin{eqnarray}
& & F^1_1 \Big(2\big(A+\frac{h'}{h}\big)(Ah+h')x^5
+\big(h''-\frac{h'}{h}-Ah'-A^2h\big)x^3
\nonumber\\
& & \label{pol-x}+\big(h''+\frac{h'}{h}+Ah'\big)x \Big).
\end{eqnarray}
where $x=\frac{z_1}{W}$.  Notice that we may suppose without loss of
generality that  $x= O(1)$. Otherwise, there exists some constant
$\alpha<1$ so that  $x\le \alpha$ what implies the estimate
\[
(1-\alpha^2)z_1^2\le \alpha^2 h^2.
\]
Thus, fixing $A=A_0$ in (\ref{pol}), the coefficients in $x$ are
uniformly bounded for $x=O(1)$. This implies that the the expression
in (\ref{pol}) is $O(\varepsilon)$ for some very small
$\varepsilon>0$. Thus, we conclude using the inequality $\psi\ge
\sum_i F^i_i a^i_i$ that (\ref{desfinal}) may be written as
\begin{equation}
\label{desfinal2} O(\varepsilon) W^2- \big(\psi_z
+\frac{h'}{h}\psi\big)z_1W^2 \le \psi_1 W^2 + A_0 h^2\psi+h'h\psi.
\end{equation}
The hypothesis $(c)$ in Theorem 1 may be stated as
\begin{equation}
\psi_z + \frac{h'}{h}\psi\le 0.
\end{equation}
Then if we choose
\[ \varepsilon \ll \frac{1}{W^2},
\]
an estimate for $W|_{\bar u}$  follows from (\ref{desfinal2}).

In both cases, by definition of the function $\chi$, a bound for
$z_1(\bar u)$ implies an uniform bound for $\nabla'z$. This
completes the proof of the Proposition  \ref{est-nabla}.

\section{Hessian bounds}

This section is devoted to the proof of Hessian estimates. We will
show that the terms of the second fundamental form are bounded by
above. Since we already have $C^0$ and $C^1$ estimates, then this
information allow us to obtain the Hessian estimates.

With this purpose in mind, we  define the following function on the
unit tangent bundle of $\Sigma$:
\begin{equation}
\tilde{\zeta}(u,\xi)=B(\xi,\xi)\, e^{\varphi(\tau)-\beta\eta},
\end{equation}
where $u\in M$, $\xi$ is an unit tangent vector to $\Sigma$ at
$(z(u),u)$, the functions $\tau$ e $\eta$ are defined in (\ref{t}),
the constant $\beta>0$ will be chosen later and the real function
$\varphi$ is defined as follows. Notice that by definition the
function  $\tau$ is bounded by constants depending on bounds for $z$
and $\nabla'z$. Hence, it is possible to choose $a>0$ so that
$\tau\geq 2a$. Thus,  we define
\begin{equation}
\varphi(\tau)=-\ln(\tau-a).
\end{equation}
Hence, one has differentiating with respect to $\tau$
\begin{equation}
\ddot
\varphi-(1+\epsilon)\dot\varphi^{2}=\frac{1}{(\tau-a)^{2}}-\frac{1
+\epsilon}{(\tau-a)^{2}}=-\frac{\epsilon}{(\tau-a)^{2}}<0
\end{equation}
and by the choice of $a$ given an arbitrary  positive constant $C$,
one has
\begin{eqnarray*}
& & -(1+\dot\varphi\tau)+C(\ddot\varphi-(1+\epsilon)\dot\varphi^{2})
\ge \hat C,
\end{eqnarray*}
for some positive constant $\hat C$ depending on bounds for $z$ and
$\nabla'z$.

We suppose that the maximum of  $\tilde{\zeta}$ is attained at a
point $\bar u$ and along the direction $\bar \xi$ tangent to $\bar
X=(z(\bar u),\bar u)$. We may choose a geodesic orthonormal
reference frame $E_a$ around $\bar X$ as defined in Section 2  so
that $\omega^k_i |_{\bar X}=0 $. One may rotate this frame in such a
way that $\bar \xi=E_1$ at $\bar X$. We then consider the local
function  $a_{11}=B(E_1,E_1)$. Thus we easily verifies that    the
function
\begin{equation}
\zeta(p)=a_{11}\, e^{\varphi(\tau)-\beta\eta}
\end{equation}
attains maximum at  $\bar X$. Thus, it holds at $\bar u$
\begin{eqnarray}\label{11}
0 =(\ln
\zeta)_i=\frac{a_{11;i}}{a_{11}}+\dot\varphi\tau_i-\beta\eta_i
\end{eqnarray}
and the Hessian matrix with components
\begin{eqnarray*}
(\ln \zeta)_{i;j}
=\frac{a_{11;ij}}{a_{11}}-\frac{a_{11;i}a_{11;j}}{a^{2}_{11}}
+\dot\varphi\tau_{i;j}+\ddot\varphi\tau_i\tau_j-\beta\eta_{i;j}
\end{eqnarray*}
is negative-definite. Thus
\begin{eqnarray}\label{12}
 F^{ij}(\ln \zeta)_{ij} &=&\frac{1}{a_{11}}F^{ij}a_{11;ij}
-\frac{1}{a^{2}_{11}}F^{ij}a_{11;i}a_{11;j}
+\dot\varphi F^{ij}\tau_{i;j}\nonumber\\
&  + &\ddot\varphi F^{ij}\tau_i\tau_j -\beta F^{ij}\eta_{ij}\le 0
\end{eqnarray}
It is clear that $a_{11}$ is the greatest eigenvalue of $B$ and
therefore $a_{1i}=0$ for $i\neq 1$. Thus, we may rotate the
orthogonal complement of $E_1$ so that in the resulting frame the
matrix $(a_{ij})$ is diagonal at $\bar X$. By Lemma \ref{deriveSm},
it results that $(F^{ij})$ is also diagonal with $F^{ii}=f_{i}$. We
denote $\lambda_{i}=a_{ii}(\bar u)$ and choose indices in such a way
that
$$\lambda_{1}\geq\lambda_{2}\geq\cdot\cdot\cdot\geq\lambda_{n}.$$
Moreover, we assume without loss of generality that $\lambda_{1}>1$
at $\bar u$. Thus, according Lemma \ref{deriveSm}, we have
$$ f_{1}\leq f_{2}\leq\cdot\cdot\cdot\leq f_{n}.$$
From (\ref{12}) one then gets
\begin{eqnarray}\label{13}
\sum_i \Big(\frac{1}{\lambda_{1}}f_{i}a_{11;ii}
-\frac{1}{\lambda_{1}^{2}}f_{i}|a_{11;i}|^{2}+\dot\varphi
f_{i}\tau_{i;i}+\ddot\varphi f_{i}|\tau_i|^{2}  -\beta
f_{i}\eta_{i;i}\Big)\le 0
\end{eqnarray}
Now, we differentiate covariantly with respect to the metric
$(g_{ij})$ in $\Sigma$ the equation  (\ref{eq:main1}) in the
direction of $E_{1}$ obtaining
$$F^{ij}a_{ij;1}=\psi_1$$
and differentiating again
\begin{equation}\label{14}
F^{ij}a_{ij;11} +F^{ij,kl}a_{ij;1}a_{kl;1}=\psi_{1;1}.
\end{equation}
From Ricci identity in Lemma 1 and using the fact that
$\delta(f)\le\sum_{i}f_{i}\lambda_{i}\le f= \psi$  we have
\begin{eqnarray*}
F^{ij}a_{ij;11}  &\le & -\lambda_{1}^{2}\delta+|\bar
R_{1010}|\psi+\sum_{i}\big(f_{i}a_{11;ii}+\lambda_{1}f_{i}\lambda_{i}^{2}
+\lambda_{1}f_{i}\bar R_{i0i0}
\\
&  + &  f_{i}\bar R_{i0i0;1}-f_{i}\bar R_{1010;i}\big).
\end{eqnarray*}
Combining this expression and  (\ref{14}) and replacing the
resulting expression in (\ref{13}) one has
\begin{eqnarray*}
& &
\frac{\psi_{1;1}}{\lambda_{1}}+\frac{1}{\lambda_{1}}\big(\delta\lambda_{1}^{2}-\psi
|\bar R_{1010}|\big) -\frac{1}{\lambda_{1}}F^{ij,kl}a_{ij;1}a_{kl;1}
-\sum_{i}f_{i}\lambda_{i}^{2}\\
& &-\sum_{i}f_{i}\bar R_{i0i0}
-\frac{1}{\lambda_{1}}\sum_{i}f_{i}\big(\bar R_{i0i0;1}-\bar
R_{1010;i}\big) +\sum_i
\Big(-\frac{1}{\lambda_{1}^{2}}f_{i}|a_{11;i}|^{2}\\
& & +\dot\varphi f_{i}\tau_{i;i}+\ddot\varphi f_{i}|\tau_i|^{2} -
\beta f_{i}\eta_{i;i}\Big)\le 0.
\end{eqnarray*}
From (\ref{ddn}) we have at  $\bar X$
\begin{eqnarray*}
\beta\sum_{i}f_{i}\eta_{i;i} = \beta\sum_i\big(\tau f_{i}a_{ii}- h'
f_{i}g_{ii}\big) \le \beta\big(\tau \psi- h' T\big),
\end{eqnarray*}
where $T=\sum_{i}f_{i}$. From  (\ref{ddt}) and denoting
\[
\bar R_{ki}: = \ae \bar R(E_k, E_i)E_i, N\ad = \bar\Omega^0_i (E_k,
E_i)
\]
and using that $\dot\varphi<0$ it holds at $\bar X$ that
\begin{eqnarray*}
\dot  \varphi\sum_{i}f_{i}\tau_{i;i} & \ge &
-\dot\varphi\Big(\sum_{i,k}\eta^k f_{i}a_{ii;k}
+\sum_{i,k}\eta^k\bar R_{ki} f_{i}\Big)-\dot\varphi\tau
\sum_{i}f_{i}\lambda_{i}^{2}\\
&+ &\dot\varphi h'\psi\\
&  = & -\dot \varphi\Big(\sum_{k}\eta^k \psi_k +\sum_{i,k}\eta^k
\bar R_{ki}f_{i}\Big)-\dot\varphi \tau
\sum_{i}f_{i}\lambda_{i}^{2}+\dot\varphi h'\psi.
\end{eqnarray*}
Using  (\ref{des3}) and estimating the ambient curvature terms by
constants $C_k$ terms one obtains from Lemma 8
\begin{eqnarray*}
-\sum_k\dot\varphi\big(\eta^k(\psi_k+C_kT)\big)+\dot \varphi h' \psi
\geq -|\dot\varphi|(C+CT).
\end{eqnarray*}
Therefore, we have
\begin{eqnarray*}
\dot \varphi \sum_{i}f_{i}\tau_{i;i} \geq
-|\dot\varphi|(C+CT)-\dot\varphi\tau \sum_{i}f_{i}\lambda_{i}^{2}.
\end{eqnarray*}
Now, we suppose without loss of generality that
\begin{equation*}
\lambda_{1}\geq\frac{1}{C}\sum_{i}|R_{i0i0;1}-R_{1010;i}|,
\end{equation*}
for some $C>0$. Moreover, supposing that $\lambda_1 \ge 1$ one has
\[
-\frac{1}{\lambda_1}\psi |\bar R_{1010}| \ge-C  \quad\textrm{  and
}\quad \frac{\psi_{1;1}}{\lambda_1}\ge -C
\]
for some positive constant $C$.  Finally one has
\[
-\sum_i f_i \bar R_{i0i0}\ge -T \max_i |\bar R_{i0i0}|\ge -CT.
\]
We then  conclude from these inequalities  that
\begin{eqnarray}\label{22}
& &
-C-CT+\delta\lambda_{1}-\frac{1}{\lambda_{1}}F^{ij,kl}a_{ij;1}a_{kl;1}
-\sum_{i}f_{i}\lambda_{i}^{2} \nonumber\\
& &- \frac{1}{\lambda_{1}^{2}}\sum_{i}f_{i}|a_{11;i}|^{2}
-|\dot\varphi|(C+CT) -\dot\varphi\tau \sum_{i}f_{i}\lambda_{i}^{2}
\nonumber\\
& & +\ddot\varphi\sum_{i}f_{i}|\tau_i|^{2}
 - \beta\big(\tau \psi- h' T\big)\le 0.
\end{eqnarray}
Finally, we also have from  (\ref{11}) for any $\epsilon>0$ the
inequality
\begin{eqnarray}
\label{ineqbt} \frac{1}{\lambda_{1}^{2}}f_{i}|a_{11;i}|^{2} \le
(1+\frac{1}{\epsilon})\beta^{2}f_{i}|\eta_i|^{2}
+(1+\epsilon)\dot\varphi^{2}f_{i}|\tau_i|^{2}.
\end{eqnarray}

\noindent Now, for proceed in our analysis, we consider two cases.

\vspace{0.3cm}

\noindent {\bf $1^{st}$ Case.}  In this case, we suppose that
$\lambda_{n}\leq-\theta\lambda_{1}$ for some positive constant
$\theta$ to be chosen later.

Replacing the sum of terms in (\ref{ineqbt}) in the inequality
(\ref{22}) and  using Lemma 8 one has after grouping terms in $T$
\begin{eqnarray*}
& & \delta\lambda_{1} -C-C|\dot\varphi|
-\frac{1}{\lambda_{1}}F^{ij,kl}a_{ij;1}a_{kl;1}
\\
& & -\big(C+C|\varphi'| -h'\beta
+C(1+\frac{1}{\epsilon})\beta^{2}\big)T\\
& & -(1+\dot\varphi\tau)\sum_{i}f_{i}\lambda_{i}^{2}
+\big(\ddot\varphi-(1+\epsilon)\dot\varphi^{2}\big)\sum_{i}f_{i}|\tau_i|^{2}-\beta
\tau \psi\le 0.
\end{eqnarray*}
Using (\ref{dt}) and the fact that  $\big(a_{ij}\big)$ is diagonal
at $\bar X$ and Lemma 8 we calculate
\begin{eqnarray}
\sum_{i}f_{i}|\tau_i|^{2}=\sum_{i}f_{i}\lambda_{i}^{2}|\eta_i|^{2}\leq
C\sum_{i}f_{i}\lambda_{i}^{2}.
\end{eqnarray}
Hence, we get
\begin{eqnarray}\label{des22}
&   &  \delta\lambda_{1}-C - C|\dot\varphi|
-\frac{1}{\lambda_{1}}F^{ij,kl}a_{ij;1}a_{kl;1}
\nonumber\\
&  & - \big(C+C|\dot\varphi|-h'\beta
+C(1+\frac{1}{\epsilon})\beta^{2}\big)T\nonumber
\\
& & +
\big(-(1+\dot\varphi\tau)+C(\ddot\varphi-(1+\epsilon)\dot\varphi^{2})\big)
\sum_{i}f_{i}\lambda_{i}^{2}-\beta\tau\psi\le 0.
\end{eqnarray}
Now, using the concavity of $F$ we may discard the third term in the
left-hand side of (\ref{des22}) since it is non-negative obtaining
\begin{eqnarray*}
-C_1(\beta) -C_{2}(\beta)T+\delta\lambda_1+\hat
C\sum_{i}f_{i}\lambda_{i}^{2}\le 0,
\end{eqnarray*}
where $C_1$ depends linearly on $\beta$ and $C_2$ depends
quadratically on $\beta$.  Since $f_{n}\geq\frac{1}{n}T$, we have
\[
\sum_{i}f_{i}\lambda_{i}^{2}\geq f_{n}\lambda_{n}^{2}\geq
\frac{1}{n}\theta^{2}T\lambda_{1}^{2}.
\]
Thus it follows that
\begin{eqnarray}
\label{polT} -C_1 - C_2 T + \delta\lambda_{1}+\hat
C\frac{1}{n}\theta^{2}T\lambda_{1}^{2}\le 0.
\end{eqnarray}
This inequality shows that  $\lambda_{1}$ has an upper bound. In
fact, if we assume without loss of generality that $\lambda_1\ge
\bar C$ for some positive constant $\bar C$,  the coefficients of
the terms in $T$ in (\ref{polT}) have a nonnegative sum. Thus,
discarding these terms, one gets
\[
\lambda_1 \le \frac{C_1}{\delta}.
\]

\vspace{0.2cm}

\noindent {\bf $2^{nd}$ Case:} In this case, we assume that
$\lambda_{n}\geq-\theta\lambda_{1}$. Hence,
$\lambda_{i}\geq-\theta\lambda_{1}$. We then group the indices in
$\{1,...,n\}$ in two sets $I_1=\{j;f_{j}\leq4f_{1}\}$ and
$I_2=\{j;f_{j}>4f_{1}\}$. Using (\ref{ineqbt}) we have for $i\in
I_1$
\begin{eqnarray*}
\frac{1}{\lambda_{1}^{2}}f_{i}|a_{11;i}|^{2}
&\leq&(1+\epsilon)\dot\varphi^{2}f_{i}|\tau_i|^{2}
+C(1+\frac{1}{\epsilon})(\beta)^{2}f_{1}.
\end{eqnarray*}
Therefore, it follows from (\ref{22}) that
\begin{eqnarray*}
& & -C-CT +\delta\lambda_{1}
-\frac{1}{\lambda_{1}}F^{ij,kl}a_{ij;1}a_{kl;1}
-\big(1+\dot\varphi\tau
\big)\sum_{i}f_{i}\lambda_{i}^{2}\\
& & - \frac{1}{\lambda_{1}^{2}}\sum_{j\in
I_2}f_{j}|a_{11;j}|^{2}-|\dot\varphi|(C+CT)
 +\big(\ddot\varphi-(1+\epsilon)\dot\varphi^{2}\big)\sum_{i}f_{i}|\tau_i|^{2}\\
& & -  C(1+\frac{1}{\epsilon})\beta^{2}f_{1}-\beta\big( \tau\psi-h'
T\big)\le 0.
\end{eqnarray*}
Notice that we had summed up to the inequality the non-positive
terms
\[
-(1+\epsilon)|\dot\varphi|^2\sum_{i\in I_2} f_i |\tau_i|^2
\]
Using Lemma 8, one has
\[
|\tau_i|^2 =|\lambda_i \eta_i|^2\le C\lambda_i^2
\]
and as we had seen above one  may prove that
\begin{equation}
\label{ineqphi}
-\big(1+\dot\varphi\tau\big)\sum_{i}f_{i}\lambda_{i}^{2}
+\big(\ddot\varphi-(1+\epsilon)\dot\varphi^{2}\big)\sum_{i}f_{i}|\tau_i|^{2}\geq
\hat C\sum_{i}f_{i}\lambda_{i}^{2}
\end{equation}
for some positive constant $\hat C$. Thus we have
\begin{eqnarray*}
& & -C-CT +\delta\lambda_{1}
-\frac{1}{\lambda_{1}}F^{ij,kl}a_{ij;1}a_{kl;1} +\hat
C\sum_{i}f_{i}\lambda_{i}^{2}
\\
& & - \frac{1}{\lambda_{1}^{2}}\sum_{j\in
I_2}f_{j}|a_{11;j}|^{2}-|\dot\varphi|(C+CT)
-C\big(1+\frac{1}{\epsilon}\big)\beta^{2}f_{1}\\
& &-\beta\big( \tau\psi- h' T\big)\le 0.
\end{eqnarray*}
Denoting $\bar R_{j1}=\bar\Omega^0_1(E_j,E_1)$ one has by Lemma
\ref{deriveSm} and the fact that $1\notin I_2$ and using Codazzi's
equation
\begin{eqnarray}
\label{FCod}
 -\frac{1}{\lambda_{1}}F^{ij,kl}a_{ij;1}a_{kl;1} \ge
-\frac{2}{\lambda_{1}}\sum_{j\in
I_2}\frac{f_{1}-f_{j}}{\lambda_{1}-\lambda_{j}}\big(a_{11;j} +\bar
R_{j1}
\big)^{2}.
\end{eqnarray}
Following \cite{QYYL}, we may verify that choosing
$\theta=\frac{1}{2}$ it holds that  for all $j\in I_2$ it holds that
\begin{equation}
\label{fj} -\frac{2}{\lambda_{1}}\frac{f_{1}-f_{j}}{\lambda_{1}
-\lambda_{j}}\geq\frac{f_{j}}{\lambda_{1}^{2}}.
\end{equation}
Considering the inequalities (\ref{FCod}) and (\ref{fj}) and using
(\ref{ineqphi}) one has
\begin{eqnarray*} & &   -C-CT+\delta\lambda_{1} +\sum_{j\in
I_2}\frac{f_{j}}{\lambda_{1}^{2}}a_{11;j}^{2} +2\sum_{j\in
I_2}\frac{f_{j}}{\lambda_{1}^{2}}a_{11;j}\bar R_{j1}
\\
& & +\hat C\sum_{i}f_{i}\lambda_{i}^{2} -\sum_{j\in
I_2}\frac{f_{j}}{\lambda^2_1}a_{11;j}^{2}-|\dot\varphi|(C+CT)
\\
& & -C(1+\frac{1}{\epsilon})\beta^{2}f_{1}-\beta\big( \tau\psi- h'
T\big)\le 0.
\end{eqnarray*}
Hence one obtains
\begin{eqnarray*}
& & -C-CT+\delta\lambda_{1} +2\sum_{j\in
I_2}\frac{f_{j}}{\lambda_{1}}(-\dot\varphi\tau_j +\beta\eta_j)\bar
R_{j1}
\\
& & +\hat C\sum_{i}f_{i}\lambda_{i}^{2}-|\dot\varphi|(C+CT)
-C(1+\frac{1}{\epsilon})\beta^{2}f_{1}\\
& &-\beta\big(\tau\psi- h' T\big)\le 0.
\end{eqnarray*}
We now estimate using that $\dot\varphi<0$ and that $\lambda_j\le
\lambda_1$ and $-\lambda_j\le \theta\lambda_1<\lambda_1$
\begin{eqnarray*}
& & 2\frac{f_{j}}{\lambda_{1}}(-\dot\varphi\tau_j)\bar R_{j1}\ge
2\frac{f_{j}}{\lambda_1}\dot\varphi|\lambda_j||\eta_j\bar R_{j1}|\ge
2f_{j}\dot\varphi|\eta_j\bar R_{j1}|.
\end{eqnarray*}
We also may suppose without loss of generality that  it holds that
\[
\lambda_1\ge \frac{3| \eta_j\bar R_{j1}|}{h'}
\]
for all $j\in I_2$. Thus, these inequalities imply that
\begin{eqnarray*}
& & -C-CT+\delta\lambda_{1} +2\sum_{j\in I_2}f_{j}\dot\varphi
|\eta_j
\bar R_{j1}| -2\frac{\beta h'}{3} T\\
& & +\hat C\sum_{i}f_{i}\lambda_{i}^{2} - |\dot\varphi|(C+CT) -
C\big(1+\frac{1}{\epsilon}\big)\beta^{2}f_{1} \\
& &- \beta\big( \tau\psi- h' T\big)\le 0.
\end{eqnarray*}
Since $\sum_{j\in I_2}f_{j}\leq T$, $|\eta_j\bar R_{j1}|\leq C$,
$\dot\varphi<0$  one has
\begin{eqnarray*}
-C- \big(C+C|\dot\varphi| + 2\beta\frac{ h'}{3}-\beta h'\big)
T-C\big(1+\frac{1}{\epsilon}\big)\beta^{2}f_{1}+\delta\lambda_{1}
+\hat C f_{1}\lambda_{1}^{2}\le 0.
\end{eqnarray*}
Choosing $\beta>0$ sufficiently large the term in $T$ is positive
and we may discard it obtaining
\begin{equation}
-C-C_2(\beta)f_1+\delta\lambda_{1}+\hat C f_{1}\lambda_{1}^{2}\le 0,
\end{equation}
where $C_2$ depends quadratically on $\beta$. Reasoning as above,
one concludes that this inequality gives an upper bound for
$\lambda_1$.

\section{The proof of the Theorem}

To prove the theorem we are going to use the degree theory for
nonlinear elliptic partial differential equations developed by Yan
Yan Li. We refer the reader to \cite{YYL}.

In Sections 3, 5 and 6 above, it is  proved that admissible $C^4$
function $z$ which solve the equation $\Upsilon(s,z)=0$ for some
$0\le s\le 1$ satisfy the following bounds
\begin{equation}
\label{C0}
 t_-< z(u)< t_+, \quad u\in M
\end{equation}
and
\begin{equation}
\label{C2} |z|_2 \le C
\end{equation}
for some positive constant $C$ which depends on $n,t_-,t_+$ and
$\psi$. Then the $C^{4,\alpha}$ estimate for some $\alpha\in[0,1]$
follows from (\ref{C2}) and from the results of  L. C. Evans e N. V.
Krylov as stated in Theorem 17.16 in \cite{GT}.  One has
\begin{equation}
\label{C4} |z|_{4,\alpha}<C
\end{equation}
for some constant $C>0$.

Fixed that  $\alpha$ we denote by $C^{4,\alpha}_{a}(M)$ the subset
of $C^{4,\alpha}(M)$ consisting of admissible functions for $F$ and
define as in Section 2 the homotopy
\begin{equation}
\Upsilon(s,\,\cdot\,):C^{4,\alpha}_{a}(M)\rightarrow
C^{2,\alpha}(M), \quad 0\le s\le 1
\end{equation}
and we consider the family of equations
$\Upsilon(s,z)= 0$.
In order to apply degree theory, we need to prove certain assertions
which are intermediate steps in the method.

It is easy to see in view of the $C^0$ and $C^1$ estimates that
there exists $\hat C >0$ for which
\begin{equation}
\hat C\le \Psi(s,z(u),u)\le \frac{1}{\hat C}, \quad u\in M,
\end{equation}
for $0\le s\le 1$ and any $z\in C^{4,\alpha}(M)$ satisfying
(\ref{C0}) and (\ref{C4}). Now, if $z\in C^{4,\alpha}_a(M)$ solves
$\Upsilon(s,z)=0$ for some $0\le s\le 1$, then
\[
F(a_{ij}(z))= \Psi(s,z(u),u)
\]
and obviously
\begin{equation}
\label{boundF} \hat C\le F(a_{ij}(z(u)))\le \frac{1}{\hat C},\quad
u\in M.
\end{equation}
However, we may verify that there exists some open bounded set
$V\subset \Gamma$ with $\bar V \subset \Gamma$ such that if
\[
\hat C\le f(\lambda_1(z(u)),\ldots, \lambda_n(z(u)))\le
\frac{1}{\hat C}
\]
then
\begin{equation}
\lambda(z(u))\in V.
\end{equation}
In particular, by (\ref{boundF}) we conclude that the matrix
$(a_{ij}(z))$ satisfies
\begin{equation}
\label{boundlambda} \lambda(a_{ij}(z))\in V.
\end{equation}
We then define the open set $\mathcal{O}$ in $C^{4,\alpha}_a(M)$
consisting of the admissible functions satisfying (\ref{C0}),
(\ref{C4}) and (\ref{boundlambda}). Thus, our reasoning above shows
that any admissible solution $z$ of $\Upsilon(s,z)=0$ for some $0\le
s\le 1$ is contained in $\mathcal{O}$. In particular, we conclude
that
\begin{equation}
\Upsilon(s,\,\cdot\,)^{-1}(0)\cap \partial \mathcal{O}=\emptyset,
\quad 0\le s\le 1.
\end{equation}
Thus, according to Definition 2.2 in \cite{YYL}  the degree
$\textrm{deg}(\Upsilon(s,\,\cdot\,),\mathcal{O},0)$ is well-defined
for all $0\leq s\leq1$.

Proposition 6 shows that $z_0=t_0$ is the unique admissible solution
to $\Upsilon(0,z)=0$ in $C^{4,\alpha}_{a}(M)$. We must prove that
the Frech\'et derivative $\Upsilon_{z}(0,z_{0})$ calculated around
$z_0$ is an invertible operator from $C^{4,\alpha}(M)$ to
$C^{2,\alpha}(M)$. One computes
\begin{eqnarray*}
\Upsilon(0,\rho z_{0}) =F(a_{ij}(\rho z_{0})) -\phi(\rho
t_{0})k(\rho t_{0})= k(\rho t_0)-\phi(\rho t_0)k(\rho t_0)
\end{eqnarray*}
and using the fact that $\phi(t_0)=1$ and that $\phi'(t_0)<0$
\begin{eqnarray*}
\Upsilon_{z}(0,z_{0})\cdot z_{0} =\frac{\dd}{\dd\rho}\Upsilon(0,\rho
z_{0})|_{\rho=1}=-\phi'(t_0)k(t_0)>0
\end{eqnarray*}
On the other hand, since obviously $\nabla' z_0 =0$ and
$\nabla'^{2}z_0=0$, then  $\Upsilon_z (0,z_0)\cdot z_0$ is just a
multiple of the zeroth order term in $\Upsilon_z(0,z_0)$.  We
conclude that $\Upsilon_{z}(0,z_{0})$ is an invertible negatively
elliptic operator.

We finally calculate
$\textrm{deg}(\Upsilon(1,\,\cdot\,),\mathcal{O},0)$. From
Proposition 2.2 in \cite{YYL}, it follows that
$\textrm{deg}(\Upsilon(s,\,\cdot \,),O,0)$ is independent  from $s$.
In particular,
\[
\textrm{deg}(\Upsilon(1,\,\cdot\,),\mathcal{O},0)=\textrm{deg}(\Upsilon(0,\,\cdot\,),\mathcal{O},0).
\]
On the other hand, we had just proved that the equation
$\Upsilon(0,z)=0$ has an unique admissible solution $z_0$ and that
the linearized operator $\Upsilon_{z}(0,z_{0})$ is invertible. Thus,
by Proposition 2.3 in \cite{YYL} one gets
\[
\textrm{deg}(\Upsilon(0,\,\cdot\,),\mathcal{O},0)=\textrm{deg}(\Upsilon_{z}(0,z_{0}),\mathcal{O},0)=\pm
1.
\]
Therefore,
\[
\textrm{deg}(\Upsilon(1,\,\cdot\,),\mathcal{O},0)\neq 0.
\]
Thus, the equation $\Upsilon(1,z)=0$ has at least one solution $z\in
O$. This completes the proof of the theorem.

\vspace{1cm}

\noindent  Francisco J. de  Andrade\\
\noindent Universidade Federal de Campina Grande \\
\noindent Centro de Forma\c c\~ao de Professores\\
\noindent Campus de Cajazeiras\\
\noindent  Cajazeiras -- Para\'\i ba\\
\noindent  58.900-000 -- Brazil

\vspace{5mm}

\noindent Jo\~ao Lucas M. Barbosa and Jorge H. S. de Lira \\
\noindent Departamento de Matem\'atica \\
\noindent Universidade Federal do Cear\'a\\
\noindent Bloco 914 -- Campus do Pici\\
\noindent Fortaleza -- Cear\'a\\
\noindent 60455-760 -- Brazil\\
\noindent joaolucasbarbosa@gmail.com\\
\noindent jorge.lira@pq.cnpq.br

\end{document}